\documentclass[10pt,a4]{article}
\usepackage{latexsym}
\usepackage{amsmath}
\usepackage{amssymb}
\usepackage{amsthm}
\usepackage{amscd}
\usepackage{mathrsfs}
\usepackage[all]{xy}

\newtheorem{theorem}{Theorem}[section]
\newtheorem{lemma}[theorem]{Lemma}
\newtheorem{corollary}[theorem]{Corollary}

\theoremstyle{definition}

\newtheorem{definition}[theorem]{Definition}

\theoremstyle{remark}
\newtheorem{notation}{Notation}

\makeatletter

\@addtoreset{equation}{section}
\makeatother

\def\La{\Lambda}

\newcommand{\cF}{{\cal F}}
\newcommand{\cK}{{\cal K}}
\newcommand{\cH}{{\cal H}}

\renewcommand{\eqref}[1]{(\ref{#1})}

\renewcommand{\bigskip}{\vspace{0.2cm}}

\begin{document}

\title{On localizations of the characteristic classes
of $\ell$-adic sheaves of rank 1}

\maketitle

\begin{center}
{\sc Takahiro Tsushima}\\
{\it
Graduate School of Mathematical Sciences,
The University of Tokyo, 3-8-1 Komaba, Meguro-ku
Tokyo 153-8914, JAPAN}
\\
E-mail: tsushima@ms.u-tokyo.ac.jp
\end{center}

\begin{abstract}
The Grothendieck-Ogg-Shafarevich formula 
is generalized to any dimensional scheme by Abbes-Kato-Saito 
in \cite{KS} and \cite{AS}. 
In this paper, we introduce two methods
of localization of the characteristic classes for sheaves of rank 1
and compare them.
As a corollary of this comparison, 
we obtain a refinement of a formula proved by Abbes-Saito in \cite{AS}
 without denominator
for a smooth sheaf of rank 1 which is clean with respect to the boundary.
\end{abstract}

\section{Introduction}
\noindent
The Grothendieck-Ogg-Shafarevich formula is a formula
calculating the Euler-Poincar\'{e} number of an $\ell$-adic sheaf
on a curve. A.~Abbes, K.~Kato and T.~Saito generalized this formula to any dimensional case.
To generalize this formula K.~Kato and T.~Saito
define the Swan class which is produced by the wild ramification
of an $\ell$-adic sheaf using logarithmic blow-up and alteration
in \cite{KS}. This invariant realizes a deep insight of S.~Bloch that the ramification of a higher dimensional arithmetic scheme produces a $0$-cycle 
class on its boundary. 
K.~Kato and T.~Saito calculated the Euler-Poincar\'{e} number
of an $\ell$-adic sheaf
on any dimensional scheme
in terms of the Swan class in \cite{KS}.
A.~Abbes and T.~Saito refined this result using
the characteristic class of an $\ell$-adic sheaf in \cite{AS}.
They compared the characteristic class with two invariants produced by the wild ramification. 
One is the Swan class mentioned above.
We call the comparison of the characteristic class with
the Swan class the Abbes-Saito formula.
The other is the 0-cycle
class $c_{\cF}$ defined by Kato in \cite{K1}
for a smooth $\La$-sheaf of {\rm rank} one which is clean with respect to the boundary 
where $\La$ is a finite commutative $\mathbb{Z}_l$-algebra.
We call the comparison  of the characteristic class with
the Kato 0-cycle class the 
Abbes-Kato-Saito formula.

We prove a localized version of the Abbes-Saito formula 
in \cite{T}
assuming the strong resolution of singularities,
which is a refinement of their formula.
We call this formula the localized Abbes-Saito formula.
We use the localized characteristic class of an $\ell$-adic sheaf
to formulate it.
The localized characteristic class
is a lifting of the characteristic class to the \'{e}tale cohomology group supported on the
boundary locus and is defined in \cite[Section 5]{AS}.
As an application of the localized Abbes-Saito formula, we proved the
Kato-Saito conductor formula in characteristic $p>0$
in \cite{T}. This was the main motivation to consider 
the localized characteristic class and the main application
of the localized Abbes-Saito formula. 
At the conference, I reported these results.

In this paper, we study two methods of localization of the characteristic classes for sheaves of rank 1. To 
refine the Kato-Saito conductor formula in characteristic $p>0$
is the main purpose to consider these localizations.
For this purpose, first we define a localization of
the characteristic class
as a cohomology class with support on the
wild locus in Section 2 using logarithmic blow-up, and we call it
the logarithmic localized characteristic class.

Recently T.~Saito introduced a notion of non-degeneration of a
$\La$-sheaf in \cite{S}.
This notion is a natural generalization of the notion of cleanness
of a $\La$-sheaf of rank 1 defined by K. Kato to higher rank.
He calculated the characteristic class of a smooth $\La$-sheaf which
is non-degenerate with respect to the boundary.

In Section 3,
we define a further localization 
of the characteristic class
as a cohomology class 
with support on the nonclean locus
for a smooth $\La$-sheaf of rank 1
inspired by an idea of T.~Saito in \cite{S}.
We call it the 
nonclean localized characteristic class of
an $\ell$-adic sheaf.
Our main theorem (Theorem \ref{theorem}) in this paper is the comparison of 
the logarithmic localized characteristic class
 with the nonclean localized
characteristic class.
As a corollary (Corollary \ref{corollary}), we obtain an equality of the logarithmic localized characteristic class of a $\La$-sheaf of rank 1
which
is clean with respect to the boundary
and the Kato 0-cycle class in the \'{e}tale cohomology group supported on the wild locus without denominator. This equality refines 
the Abbes-Kato-Saito formula mentioned above.

I would like to thank Professor T.~Saito for introducing this subject
to me, suggesting that there exists a cohomology class
with support on the wild locus which refines the characteristic class, encouragements and so many advices
and thank the referee for many comments on an earlier version of this paper.
\begin{notation}In this paper,
$k$ denotes a field.
Schemes over $k$ are assumed to be separated and of finite type.
For a divisor with simple normal crossings
of a smooth scheme over $k$,
we assume that the irreducible components and their intersections are also
smooth over $k.$
The letter $l$ denotes a prime number invertible in $k$
and $\La$
denotes a finite commutative $\mathbb{Z}_l$-algebra.
For a scheme $X$ over $k$,
Let $\cK_X$ denote $Rf^!\La$
where $f:X \longrightarrow {\rm Spec} k$ is the structure
map. If $f$ is smooth, the canonical class map $\La(d)[2d] \longrightarrow
\cK_X$ is an isomorphism by \cite{AS} (1.8). 
When  we say a scheme $X$ is of dimension $d$, we understand
that every irreducible component of $X$ is of dimension $d$.
\end{notation}
\section{Logarithmic localized characteristic class of a smooth $\La$-sheaf of rank 1}
In \cite{T}, we define a localization of the characteristic class
of a smooth $\La$-sheaf 
of any rank as a cohomology class with support on the wild locus which we call the logarithmic localized characteristic class.
In this section, we introduce a more elementary definition
of the logarithmic localized characteristic class of a
smooth $\La$-sheaf of rank 1.

Let $X$ be a smooth scheme of dimension $d$ over $k,$ 
$U \subset X$ an open subscheme, 
the complement $X \setminus U=\bigcup_{i \in I}D_i$
a divisor with simple normal crossings and
$j:U \longrightarrow X$ the open immersion.
We recall the definitions of the logarithmic blow-up and
the logarithmic product from
\cite[Section 2.2]{AS} and \cite[Section 1.1]{KS}.
For $i \in I$, let $(X \times X)'_i \longrightarrow X \times X$
be the blow-up at $D_i \times
 D_i$ and let
$(X \times X)_{i} \widetilde{} \subset (X \times X)'_i$
be the complement of the proper
transforms of
$D_i \times X$ and $X \times D_i.$
We define the log blow-up with respect to divisors $\{D_i\}_{i \in I}$
$$(X \times X)' \longrightarrow X \times X$$
to be the fiber product
$\coprod_{i\in I} (X \times X)'_i \longrightarrow
X \times X$ of $(X \times X)'_i (i \in I)$ over $X \times X.$
We define the log product with respect to divisors $\{D_i\}_{i \in I}$
$$(X \times X) \widetilde{} \subset (X \times X)'$$
to be the fiber product
$\coprod_{i\in I} (X \times X) \widetilde{}_i \longrightarrow
X \times X$ of $(X \times X) \widetilde{}_i (i \in I)$ over $X \times X.$
The log blow-up $(X \times X)'$ is a smooth scheme over $k$
of dimension $2d.$
\begin{lemma}\label{geom}
Let the notation be as above.
We consider the following cartesian diagram
\[\xymatrix{
\widetilde{X} \ar[r]^{\!\!\!\!\!\!\!\!\!\!\!i'}\ar[d] &
(X \times X) \widetilde{} \ar[d] \\
X \ar[r]^{\!\!\!\!\!\!\!\!\delta} &
X \times X
}
\]where the right vertical arrow 
$(X \times X) \widetilde{} \longrightarrow X \times X$
is the projection and $\delta:X \longrightarrow X \times X$
is the diagonal closed immersion.
Then the inverse image $\widetilde{X}$ is the union of 
$(\mathbb{G}_m)^{\sharp{J}}$-bundles $\{D_{J} \times (\mathbb{G}_m)^{\sharp{J}}\}_{J \subset I}$ of dimension $d$
where
$D_J$
is the intersection of $\{D_i\}_{i \in J}$ in $X.$
\end{lemma}
\begin{proof}
Let $\widetilde{X}_i$
denote the inverse image of the diagonal $X$
by the projection $(X \times X)_i \widetilde{}
\longrightarrow X \times X$ for $i \in I.$
By the definition of the log product,
$\widetilde{X}_i$ is the union of the diagonal 
$X \subset (X \times X)_i \widetilde{}$
and ${\mathbb{G}_m}_{D_i}$.
Since $(X \times X) \widetilde{}$
is the fiber product of $(X \times X)_i \widetilde{}$
($i \in I$)
over $X \times X$, $\widetilde{X}$
is the fiber product of the schemes $\{\widetilde{X}_i\}_{i \in I}$
($i \in I$) over $X.$
Therefore the inverse image $\widetilde{X}$
is the union of $(\mathbb{G}_m)^{\sharp{J}}$-bundles $\{D_{J} \times (\mathbb{G}_m)^{\sharp{J}}\}_{J \subset I}$ of dimension $d$
where
$D_J$
is the intersection of $\{D_i\}_{i \in J}$ in $X.$
Hence the assertion follows.
\end{proof}

We keep the above notation. 
For a subset $I^+ \subset I,$
we put $D^+=\bigcup_{i \in I^+}D_i \subset D$.
Let $V \subset X$ denote the complement
of $D^+$ in $X.$
Let $U'$ be the complement of $D':=\bigcup_{i \in I'}D_i$
in $X$ where $I'=I \setminus I^+.$
We have $U=U' \cap V.$
Let $j':U \longrightarrow V$ and
$j_{U'}:U' \longrightarrow X$
be the open immersions.
Let 
$(X \times X) \widetilde{}_{U'} \subset (X \times X)'_{U'}$
denote the log product and log blow-up with respect to
$\{D_i\}_{i \in I'}$ respectively.

We consider the cartesian diagram
\[\xymatrix{
(X \times X) \widetilde{} _{U'} \ar[d]_{f} &
(V \times X) \widetilde{} \ar[l]_{\widetilde{j_1}\!\!\!\!\!\!\!\!\!}\ar[d]^{f_1} &
(V \times V) \widetilde{} \ar[l]_{\widetilde{k}_1\!\!\!\!\!}\ar[d]^{\widetilde{f}} \\
X \times X &
V \times X \ar[l]_{j_1\!\!\!\!\!\!} &
V \times V \ar[l]_{k_1\!\!\!\!\!\!}\\
}
\]
where the horizontal arrows are the open immersions
and the vertical arrows are the projections.
Let 
$\delta_U:U \longrightarrow U \times U$ be the diagonal closed immersion,
and
$\widetilde{\delta}:X \longrightarrow (X \times X) \widetilde{} _{U'}$
and
$\widetilde{\delta}_{V}:V \longrightarrow (V \times V) \widetilde{}$
be the logarithmic diagonal closed immersions
induced by the universality of blow-up.
We consider the cartesian diagrams
\[\xymatrix{
\widetilde{X}_{U'} \ar[r]^{\!\!\!\!\!\!\!\!\!\!\!\!\widetilde{i}}\ar[d] &
(X \times X) \widetilde{}_{U'} \ar[d]^{f} &
(X \times X) \widetilde{}_{U'} \setminus \widetilde{X}_{U'} \ar[l]_{\!\!\!\!\!\!\!\!\!\widetilde{g}}\ar[d]\\
X  \ar[r]^{\!\!\!\!\!\!\!\!\!\delta_X} &
X \times X &
X \times X \setminus \delta_X(X) \ar[l],\\}
\]
and
\[\xymatrix{
\widetilde{V} \ar[r]^{\!\!\!\!\!\!\!\!\widetilde{i}_V}\ar[d] &
(V \times V) \widetilde{} \ar[d]^{\widetilde{f}} &
(V \times V) \widetilde{} \setminus \widetilde{V} \ar[l]_{\widetilde{g}_V}\ar[d]\\
V  \ar[r]^{\!\!\!\!\delta_V} &
V \times V &
V \times V \setminus \delta_V(V) \ar[l]\\}
\]
where $\widetilde{g}:(X \times X) \widetilde{}_{U'} \setminus \widetilde{X}_{U'} \longrightarrow
 (X \times X) \widetilde{}_{U'}$ and $\widetilde{g}_V:(V \times V) \widetilde{} \setminus \widetilde{V} \longrightarrow (V \times V) \widetilde{}$
are the open immersions.

Let $\cF$
be a smooth
$\La$-sheaf of rank 1 on $U$ which is tamely ramified along $V \setminus U$. We put $\cH_0:=\mathcal{H}om({\rm pr}_2^ \ast \cF,{\rm pr}_1^\ast \cF)$ and
 $\bar{\cH}:=R\mathcal{H}om({\rm pr}_2 ^\ast j_!\cF,R{\rm pr}_1 ^! j_!\cF)$
on $U \times U$ and $X \times X$ respectively.
Further, we put 
$\bar{\cH}_V:=R\mathcal{H}om({\rm pr}_2^\ast j'_!\cF,R{\rm pr}_1^!j'_!\cF)$
on 
$V \times V$ 
and
$\widetilde{\cH}:=({\widetilde{j}}_\ast \cH_0)(d)[2d]$ on
$(V \times V) \widetilde{}$ respectively.
There exists a unique
map 
\begin{align}\label{pull}
{\widetilde{f}}^\ast \bar{\mathcal{H}}_V
\longrightarrow
\widetilde{\cH}
\end{align}
inducing the canonical isomorphism
$R\mathcal{H}om({\rm pr}_2^\ast \cF,R{\rm pr}_1^!\cF)
\longrightarrow \cH_0(d)[2d]$ on $U \times U$
defined in \cite[Proposition 3.1.1.1]{S}.
We put
$\widetilde{\cH}_{U'}:={\widetilde{j_1}}_!
{R\widetilde{k}_1}_\ast \widetilde{\cH}$
on $(X \times X) \widetilde{} _{U'}.$
We define a map
\begin{equation}\label{11}
f ^\ast \bar{\mathcal{H}}
\longrightarrow
\widetilde{\cH}_{U'}\end{equation}
to be the composition of the following maps
$$f ^\ast \bar{\mathcal{H}} \simeq 
f ^\ast {j_1}_!R{k_1}_\ast \bar{\cH}_V
\simeq
{\widetilde{j_1}}_!f_1^\ast R{k_1}_\ast \bar{\mathcal{H}}_V
\longrightarrow
{\widetilde{j_1}}_!{R\widetilde{k}_1}_\ast {\widetilde{f}}^\ast \bar{\mathcal{H}}_V
\longrightarrow
{\widetilde{j_1}}_!{R\widetilde{k}_1}_\ast \widetilde{\cH}=\widetilde{\cH}_{U'}$$
where 
the first isomorphism is induced by the Kunneth formula,
the second and third maps are induced by the base change maps
$f ^\ast {j_1}_! \longrightarrow {\widetilde{j_1}}_!f_1^\ast $
and $f_1^\ast R{k_1}_\ast 
\longrightarrow {R\widetilde{k}_1}_\ast {\widetilde{f}}^\ast $
respectively,
 and the fourth map is induced by applying the functor
${\widetilde{j_1}}_!{R\widetilde{k}_1}_\ast $
to the map (\ref{pull}).

By the definition of $\widetilde{\cH}_{U'}$, we have an
isomorphism
$\widetilde{\delta}^\ast {\widetilde{\cH}}_{U'}
\simeq
{j_V}_! {\widetilde{\delta}_{V}}^\ast (\widetilde{j}_\ast \cH_0)(d)[2d].$
The base change map ${\widetilde{\delta}_{V}}^\ast \widetilde{j}_\ast \cH_0
\longrightarrow {j'}_\ast \delta_U^\ast \cH_0
\simeq j'_\ast \mathcal{E}nd(\cF)$
 and the trace map $\mathcal{E}nd(\cF) \longrightarrow \La_U$ induce a trace map
\begin{equation}\label{Tr}
\widetilde{\rm Tr}:\widetilde{\delta}^\ast  {\widetilde{\cH}}_{U'}
\simeq{j_V}_!{\widetilde{\delta}_{V}}^\ast (\widetilde{j}_\ast \cH_0)(d)[2d]
\longrightarrow {j_V}_!\La_V(d)[2d]={j_V}_! \cK_V
\longrightarrow \cK_X.
\end{equation}

The map (\ref{11}) induces the pull-back
\begin{equation}\label{01}
f^\ast:H_X^0(X \times X,\bar{\cH}) \longrightarrow H_{\widetilde{X}_{U'}}^0((X \times X) \widetilde{}_{U'},\widetilde{\cH}_{U'}).
\end{equation}
The canonical map $\La \longrightarrow R{\widetilde{g}}_\ast \La$
induces a map
\begin{equation}\label{13}
H_{\widetilde{X}_{U'}}^0((X \times X) \widetilde{}_{U'},\widetilde{\cH}_{U'})
\longrightarrow
H_{\widetilde{X}_{U'}}^0((X \times X) \widetilde{}_{U'},\widetilde{\cH}_{U'} \otimes R{\widetilde{g}}_\ast \La).
\end{equation}
\begin{lemma}\label{lemma1}
The canonical map
\begin{equation}\label{13'}
H_{\widetilde{X}_{U'} \setminus \widetilde{V}}^0((X \times X) \widetilde{}_{U'},\widetilde{\cH}_{U'} \otimes R{\widetilde{g}}_\ast \La) 
\longrightarrow
H_{\widetilde{X}_{U'}}^0((X \times X) \widetilde{}_{U'},\widetilde{\cH}_{U'} \otimes R{\widetilde{g}}_\ast \La)
\end{equation}
is an isomorphism.
\end{lemma}
\begin{proof}
By the localization sequence, it suffices to prove
$H_{\widetilde{V}}^i((V \times V) \widetilde{},\widetilde{\cH} \otimes {R{\widetilde{g}_V}}_\ast \La)=0$ for all $i$.
Since $\cF$ is a smooth sheaf of rank 1 which is tamely ramified along $V \setminus U$, $\widetilde{\cH}$ is a smooth $\La$-sheaf of rank 1 on $(V \times V) \widetilde{}$
by \cite[Proposition 4.2.2.1]{AS}.
Therefore the canonical map $\widetilde{\cH} \otimes {R{\widetilde{g}_V}}_\ast \La \longrightarrow {R{\widetilde{g}_V}}_\ast {{\widetilde{g}_V}}^\ast \widetilde{\cH}$ is an isomorphism by the projection formula.
By this isomorphism, we obtain
isomorphisms
$H_{\widetilde{V}}^i((V \times V) \widetilde{},\widetilde{\cH} \otimes {R{\widetilde{g}_V}}_\ast \La) \simeq H^i(\widetilde{V},R{\widetilde{i}_V}^!(\widetilde{\cH} \otimes {R{\widetilde{g}_V}}_\ast \La)) \simeq H^i(\widetilde{V},R{\widetilde{i}_V}^!{R{\widetilde{g}_V}}_\ast {{\widetilde{g}_V}}^\ast \widetilde{\cH})=0$ for all $i$. Hence the assertion follows.
\end{proof}
We consider the pull-back by $\widetilde{\delta}$
\begin{equation}\label{14}
\widetilde{\delta} ^\ast:H_{\widetilde{X}_{U'} \setminus \widetilde{V}}^0((X \times X) \widetilde{}_{U'},\widetilde{\cH}_{U'} \otimes R{\widetilde{g}}_\ast \La) \longrightarrow
H_{D^+}^0(X,\widetilde{\delta} ^\ast \widetilde{\cH}_{U'} \otimes \widetilde{\delta} ^\ast R{\widetilde{g}}_\ast \La).
\end{equation}
The trace map (\ref{Tr}) induces a map
\begin{equation}\label{15}
{\widetilde{\rm Tr}}:H_{D^+}^0(X,\widetilde{\delta} ^\ast \widetilde{\cH}_{U'} \otimes \widetilde{\delta} ^\ast R{\widetilde{g}}_\ast \La)
\longrightarrow
H_{D^+}^0(X,\cK_X \otimes \widetilde{\delta} ^\ast R{\widetilde{g}}_\ast \La).
\end{equation}

\begin{lemma}\label{lemma2}
The canonical map $\La \longrightarrow
\widetilde{\delta} ^\ast R{\widetilde{g}}_\ast \La$
induces an isomorphism
\begin{equation}\label{16}
H_{D^+}^0(X,\cK_X) \simeq
H_{D^+}^0(X,\cK_X \otimes \widetilde{\delta} ^\ast R{\widetilde{g}}_\ast \La).
\end{equation}
\end{lemma}
\begin{proof}
We consider the distinguished triangle 
$$R{\widetilde{i}}^! \La|_X(d)[2d]
\longrightarrow \cK_X \longrightarrow \cK_X \otimes
{\widetilde{\delta}}^\ast R{\widetilde{g}}_\ast \La \longrightarrow$$
by the isomorphism $\La(d)[2d]
\simeq \cK_X.$
By this triangle, it
is sufficient to prove that $H_{D^+}^{j}(X,R{\widetilde{i}}^! \La|_X)$ is zero
for $j<2d+2$. 
For a subset $J \subset I',$ let $\widetilde{X}_J=D_J \times (\mathbb{G}_m)^{\sharp{J}}$ 
denote the $(\mathbb{G}_m)^{\sharp{J}}$-bundles over $D_J.$
By Lemma \ref{geom} and the purity theorem,
we obtain isomorphisms 
$R^j{\widetilde{i}}^!\La=0$
for $j<2d$,
$R^{2d}{\widetilde{i}}^!\La \simeq \bigoplus_{{J} \subset I'}\La_{\widetilde{X}_J}(-d)$
and $R^{2d+1}{\widetilde{i}}^!\La=0.$
By these isomorphisms and $D_J=\widetilde{\delta}(X) \cap \widetilde{X}_J,$ we 
acquire
$H_{D^+}^j(X,R{\widetilde{i}}^! \La|_X)=0$
for $j<2d,$
$H_{D^+}^{2d+i}(X,R{\widetilde{i}}^! \La|_X) \simeq 
\bigoplus_{{J} \\\subset I'}H_{D^+ \cap D_J}^i(D_J, \La_{D_J}(-d))$
for $i=0,1.$
Hence the assertion follows from the purity theorem.
\end{proof}
\begin{definition}\label{Deff}
We call the image
of the element
${\rm id}_{j_!\cF} \in {\rm End}_X(j_!\cF) \simeq H_X^0(X \times X,\bar{\cH})$
under the composite of the maps
(\ref{01})-(\ref{16})
 {\it the logarithmic localized characteristic class of $j_!\cF$}
and we denote it by
$C_{D^+}^{{\rm log},0}(j_!\cF) \in H_{D^+}^0(X,\cK_X).$
\end{definition}

\section{Nonclean localization of the characteristic class
and comparison with the logarithmic localized characteristic class}
In this section, we will define a further localization
of the characteristic class as a cohomology class with support on the nonclean locus which we call the nonclean localized characteristic class.
 
Let the notation be as in the previous section.
Further, we assume that $k$ is a perfect field.
Let $R=\Sigma_{i \in I}r_iD_i$ be the Swan divisor of $\cF$ 
with respect to the boundary $D$
defined by K.~Kato in \cite{K1} and \cite{K2}.
We assume that the support of the divisor $R$ is $D^+.$
We regard $R \subset X$ as a closed subscheme
of $(X \times X)'$ by the log diagonal
map $X \longrightarrow (X \times X)'$
and 
let $\pi:(X \times X)^{[R]} \longrightarrow (X \times X)'$ denote the blow-up at $R.$
Let $\Delta_i \subset (X \times X)'$ be the exceptional divisor above $D_i \times D_i$ for each $i \in I.$ Let $(X \times X)^{(R)}$ be the complement  of the union of the proper transforms of $\Delta_i$
for $i \in I^+$ in
$\pi^{-1}((X \times X) \widetilde{}).$
We consider the cartesian diagram
\[\xymatrix{
X^{(R)} \ar[r]^{\!\!\!\!\!\!\!i^{(R)}}\ar[d] & (X \times X)^{(R)}\ar[d]^{{f^{(R)}}} \\
X \ar[r]^{\!\!\!\!\!\!\!\!\!\!\delta_X} & X \times X}
\]
where ${f^{(R)}}:(X \times X)^{(R)} \longrightarrow X \times X$
is the projection.
Let $j^{(R)}: U \times U \longrightarrow (X \times X)^{(R)}$ be the open immersion and $\delta^{(R)}:X \longrightarrow (X \times X)^{(R)}$ the extended diagonal closed immersion induced by the universality of blow-up.(See \cite[Section 4.2]{AS} and \cite[Section 2.3]{S}.)

Let ${E^{+}}$ be the complement of $(V \times V) \widetilde{}$
in $(X \times X)^{(R)}$ which is a vector bundle over $D^{+}.$
Let $T \subset D^+$ be the nonclean locus of $\cF$
which is a closed subscheme in $X,$ and 
$W$ the complement of $T$ in $X$.
(See \cite{K2} for the notion of cleanness.)
Then, the sheaf $\cF$ is clean with respect to
the boundary $W \setminus U$.
Let $R_W$ denote the restriction of $R$ to the open subscheme $W \subset X$.
We have the open subschemes $U \subset V \subset W \subset X.$
We consider the cartesian diagram
\[\xymatrix{
(X \times X)^{(R)} \ar[d]^{h}&
(V \times V) \widetilde{}  \ar[l]_{{\widetilde{j}}^{(R)}\!\!\!\!\!\!\!\!\!\!\!\!\!\!}\ar[d]^{\widetilde{k}_1} \\
(X \times X) \widetilde{}_{U'}  &
(V \times X) \widetilde{} \ar[l]_{\widetilde{j_1}\!\!\!\!\!\!} 
}
\]
where $h:(X \times X)^{(R)} \longrightarrow
(X \times X) \widetilde{}_{U'}$ is the projection and all the arrows except $h$
are the open immersions.

We put 
$\cH^{(R)}:=j^{(R)}_\ast \cH_0(d)[2d]$ on $(X \times X)^{(R)}.$
The map of functors ${{\widetilde{j}}^{(R)}}_! \longrightarrow {{\widetilde{j}}^{(R)}}_\ast$
induces a map
\begin{equation}\label{0}
h^\ast {\widetilde{\cH}}_{U'}=
h^\ast {\widetilde{j_1}}_!{R\widetilde{k}_1}_\ast {\widetilde{\cH}}
\simeq {{\widetilde{j}}^{(R)}}_! {\widetilde{\cH}} \longrightarrow {\widetilde{j}^{(R)}}_\ast {\widetilde{j}}_\ast \cH_0(d)[2d] \simeq \cH^{(R)}.
\end{equation}
By the maps (\ref{11}) and (\ref{0}), we obtain the following map
\begin{equation}\label{2}
{{f^{(R)}}}^\ast \bar{\cH} \longrightarrow \cH^{(R)}
\end{equation}
defined in \cite[Corollary 3.1.2.2]{S}.
This map induces the pull-back 
\begin{equation}\label{p-ll}
{{f^{(R)}}}^\ast :H_X^0(X \times X, \bar{\cH})
\longrightarrow
H_{X^{(R)}}^0((X \times X)^{(R)},\cH^{(R)}).
\end{equation}

There exists a unique section $e \in \Gamma(X,{{\delta}^{(R)}}^\ast j^{(R)}_\ast  \cH_0)$
lifting the identity
${\rm id}_{\cF} \in {\rm End}_U(\cF) \simeq \Gamma(U,\delta_U^\ast \cH_0).$ 
(See \cite[Definition 2.3.1]{S}.)
We consider the natural cup pairing
\begin{equation}\label{pairing1}
\cup:\Gamma(X,{\delta^{(R)}}^\ast j^{(R)}_\ast \cH_0) \times H_X^{2d}((X \times X)^{(R)},\La(d)) \longrightarrow H_X^0((X \times X)^{(R)},\cH^{(R)}).
\end{equation}
We write
$[X]$ for the image of the cycle class $[X] \in CH_d(X)$
under the cycle class map $CH_d(X)
\longrightarrow H_{X}^{2d}((X \times X)^{(R)},\La(d)).$
By the pairing (\ref{pairing1})
and the pull-back (\ref{p-ll}), we obtain an element
${f^{(R)}}^\ast {\rm id}_{j_!\cF}-e \cup [X] \in
H_{X^{(R)}}^0((X \times X)^{(R)},\cH^{(R)}).$

We have the localization sequence
$$H_{W^{(R_W)}}^{-1}((W \times W)^{(R_W)},\cH^{(R_W)}_W) \longrightarrow
H_{T^{(R)}}^0((X \times X)^{(R)},\cH^{(R)}) \longrightarrow
H_{X^{(R)}}^0((X \times X)^{(R)},\cH^{(R)}) \longrightarrow $$
\begin{equation}\label{loc}
H_{W^{(R_W)}}^0((W \times W)^{(R_W)},\cH^{(R_W)}_W) \longrightarrow
H_{T^{(R)}}^1((X \times X)^{(R)},\cH^{(R)}) \longrightarrow.....
\end{equation} 
where $T^{(R)}$ is the complement of $W^{(R_W)}$
in $X^{(R)}$ and $\cH^{(R_W)}_W$ is the restriction of ${\cH}^{(R)}$
to the open subscheme $(W \times W)^{(R_W)} \subset {(X \times X)^{(R)}}.$

\begin{lemma}\label{l1}
The canonical map
$H_{T^{(R)}}^0((X \times X)^{(R)},\cH^{(R)})
\longrightarrow
H_{X^{(R)}}^0((X \times X)^{(R)},\cH^{(R)})$ is injective.
\end{lemma}
\begin{proof}
It suffices to prove
$H_{W^{(R_W)}}^{-1}((W \times W)^{(R_W)},\cH^{(R_W)}_W)=0$
by the localization sequence (\ref{loc}).
We consider the exact sequence
$$
H_{W^{(R_W)} \setminus \widetilde{V}}^{-1}((W \times W)^{(R_W)},\cH^{(R_W)}_W) 
\longrightarrow
H_{W^{(R)}}^{-1}((W \times W)^{(R_W)},\cH^{(R_W)}_W) \longrightarrow
H_{\widetilde{V}}^{-1}((V \times V) \widetilde{},\widetilde{\cH}).$$ 
Since we have 
$H_{W^{(R_W)} \setminus \widetilde{V}}^{-1}((W \times W)^{(R_W)},\cH^{(R_W)}_W)=0$ 
by \cite[the proof of Theorem 3.1.3]{S}, it suffices to prove 
the vanishing
$H_{\widetilde{V}}^{-1}((V \times V) \widetilde{},\widetilde{\cH})=0$ by the exact sequence.
Since $\widetilde{\cH}$ is a smooth sheaf on $(V \times V) \widetilde{}$
by \cite[Proposition 4.2.2.1]{AS}, the canonical map 
${\widetilde{i}_V}^\ast \widetilde{\cH} \otimes R{\widetilde{i}_V}^! \La
\longrightarrow
R{\widetilde{i}_V}^!\widetilde{\cH}
$ is an isomorphism
by the projection formula.
By this isomorphism, we acquire isomorphisms $H_{\widetilde{V}}^{-1}((V \times V) \widetilde{},\widetilde{\cH})
\simeq H^{-1}(\widetilde{V},R{\widetilde{i}_V}^!\widetilde{\cH})
\simeq H^{-1}(\widetilde{V},{\widetilde{i}_V}^\ast \widetilde{\cH} \otimes R{\widetilde{i}_V}^! \La) \simeq H^{2d-1}(\widetilde{V},{\widetilde{i}_V}^\ast {\widetilde{j}}_\ast \widetilde{\cH}_0(d) \otimes R{\widetilde{i}_V}^! \La).$
Thereby it is sufficient to prove that
$H^{2d-1}(\widetilde{V},{\widetilde{i}_V}^\ast {\widetilde{j}}_\ast \widetilde{\cH}_0(d) \otimes R{\widetilde{i}_V}^! \La)=0.$
By Lemma \ref{geom} and the purity theorem, the sheaf 
$R^j{\widetilde{i}_V}^!\La$ is zero
for $j<2d$.
Hence the assertion follows.
\end{proof}

\begin{corollary}\label{C1}
Let the notation be as above.
Then there exists a unique element $({f^{(R)}}^\ast {\rm id}_{j_!\cF})'$ in $H_{T^{(R)}}^0((X \times X)^{(R)},\cH^{(R)})$
which goes to the element ${f^{(R)}}^\ast {\rm id}_{j_!\cF}-e \cup [X] \in H_{X^{(R)}}^0((X \times X)^{(R)},\cH^{(R)})$
by the canonical map $H_{T^{(R)}}^0((X \times X)^{(R)},\cH^{(R)})
\longrightarrow
H_{X^{(R)}}^0((X \times X)^{(R)},\cH^{(R)}).$
\end{corollary}
\begin{proof}
Since $\cF$ is clean with respect to the boundary $W \setminus U,$
the element ${f^{(R)}}^\ast {\rm id}_{j_!\cF}-e \cup [X]$
goes to zero under the restriction map 
$H_{X^{(R)}}^0((X \times X)^{(R)},\cH^{(R)})
\longrightarrow
H_{W^{(R_W)}}^0((W \times W)^{(R_W)},\cH^{(R_W)}_W)$
by \cite[Corollary 3.1.2 and Theorem 3.1.3]{S}. Therefore
the assertion follows from Lemma \ref{l1}
and the sequence (\ref{loc}).
\end{proof}
The base change map ${\delta^{(R)}}^\ast j^{(R)}_\ast  \longrightarrow j_\ast \delta_U ^\ast$ and
the trace map $\mathcal{E}nd(\cF) \longrightarrow \La_U$
induce a trace map
\begin{align}\label{Tr'}
{\rm Tr}^{(R)}:{\delta^{(R)}}^\ast \cH^{(R)} \longrightarrow j_\ast \delta_U ^\ast \cH_0(d)[2d]\simeq j_\ast \mathcal{E}nd(\cF)(d)[2d] \longrightarrow \La_X(d)[2d]\simeq \cK_X.
\end{align}

\begin{definition}\label{deg def}
We consider the following map 
\[\xymatrix{
({f^{(R)}}^\ast {\rm id}_{j_!\cF})'
\in H_{T^{(R)}}^0((X \times X)^{(R)},\cH^{(R)}) \ar[r]^{{{\delta}^{(R)}}^ \ast \!\!\!\!\!\!\!\!\!\!\!\!\!\!\!\!\!\!\!\!\!\!\!\!\!\!\!\!\!\!\!\!\!\!\!\!\!\!\!\!\!} &
H_{T}^0(X,{{\delta}^{(R)}}^ \ast \cH^{(R)}) \ar[r]^{{\rm Tr}^{(R)}\!\!\!\!\!\!\!\!\!\!\!\!\!} &
H_{T}^0(X,\cK_X)}
\] induced by the pull-back by ${\delta}^{(R)}$ and the trace map (\ref{Tr'}).
We call the  image of the element $({f^{(R)}}^\ast {\rm id}_{j_!\cF})'$
in Corollary \ref{C1}
by this composition
{\it the nonclean localized characteristic class of $j_!\cF$} and we denote it by $C_T^{\rm ncl,0}(j_!\cF).$
\end{definition}

We will compare $C_T^{\rm ncl,0}(j_!\cF)$
with $C_{D^+}^{\rm log,0}(j_!\cF)$
defined in Definition \ref{Deff} in $H_{D^+}^0(X,\cK_X).$
(Theorem \ref{theorem}.)

\begin{lemma}\label{l3}
The canonical maps
\begin{equation}\label{map1}
H_{E^+}^0((X \times X)^{(R)},\cH^{(R)} \otimes h^\ast R{\widetilde{g}}_\ast \La)
\longrightarrow
H_{X^{(R)}}^0((X \times X)^{(R)},\cH^{(R)} \otimes h^\ast R{\widetilde{g}}_\ast \La),\end{equation}
\begin{equation}\label{map3}
H_{D^+}^0((X \times X)^{(R)},\cH^{(R)} \otimes h^\ast R{\widetilde{g}}_\ast \La) \longrightarrow
H_X^0((X \times X)^{(R)},\cH^{(R)} \otimes h^\ast R{\widetilde{g}}_\ast \La),
\end{equation}
\begin{equation}\label{map4}
H_{E^+}^0((X \times X)^{(R)},h^\ast R{\widetilde{g}}_\ast \La)
\longrightarrow
H_{X^{(R)}}^0((X \times X)^{(R)},h^\ast R{\widetilde{g}}_\ast \La)
\end{equation} and
\begin{equation}\label{map2} 
H_{D^+}^{2d}((X \times X)^{(R)},h^\ast R{\widetilde{g}}_\ast \La(d))
\longrightarrow
H_{X}^{2d}((X \times X)^{(R)},h^\ast R{\widetilde{g}}_\ast \La(d))
\end{equation}
are isomorphisms.
\end{lemma}
\begin{proof}
We prove the assertions in the same way as Lemma \ref{lemma1}.
\end{proof}
We consider the composite
{\small$$H_{X^{(R)}}^0((X \times X)^{(R)},\cH^{(R)}) \longrightarrow
H_{X^{(R)}}^0((X \times X)^{(R)},\cH^{(R)} \otimes h^\ast R{\widetilde{g}}_\ast \La)\simeq
H_{E^+}^0((X \times X)^{(R)},\cH^{(R)} \otimes h^\ast R{\widetilde{g}}_\ast \La)$$}
where the first map is induced by the canonical map
$\La \longrightarrow h^\ast R{\widetilde{g}}_\ast \La$
and the second isomorphism is (\ref{map1}).
Let $({f^{(R)}}^\ast {\rm id}_{j_!\cF})^{\rm log} \in 
H_{E^+}^0((X \times X)^{(R)},\cH^{(R)} \otimes h^\ast R{\widetilde{g}}_\ast \La)$
denote the image of the element
${f^{(R)}}^\ast {\rm id}_{j_!\cF}
\in H_{X^{(R)}}^0((X \times X)^{(R)},\cH^{(R)})$
by this map.

The pull-back by ${\delta^{(R)}},$
the trace map (\ref{Tr'}) and the isomorphism (\ref{16}) induce a map
$${\rm Tr}^{(R)} \cdot {\delta^{(R)}}^\ast:H_{E^+}^0((X \times X)^{(R)},\cH^{(R)} \otimes h^\ast R{\widetilde{g}}_\ast \La)
\longrightarrow
H_{D^+}^0(X,\cK_X \otimes {\widetilde{\delta}}^\ast R{\widetilde{g}}_\ast \La) \simeq
H_{D^+}^0(X,\cK_X).$$
We write
 ${\rm Tr}^{(R)} {\delta^{(R)}}^\ast ({f^{(R)}}^\ast {\rm id}_{j_!\cF})^{\rm log} \in H_{D^+}^0(X,\cK_X)$
for the image of the element
$({f^{(R)}}^\ast {\rm id}_{j_!\cF})^{\rm log}
\in H_{E^+}^0((X \times X)^{(R)},\cH^{(R)} \otimes h^\ast R{\widetilde{g}}_\ast \La)$
under the 
map ${\rm Tr}^{(R)} \cdot {\delta^{(R)}}^\ast.$

\begin{lemma}\label{l3'}
Let the notation be as above.
We have an equality
$$C_{D^+}^{\rm log,0}(j_!\cF)={\rm Tr}^{(R)} {\delta^{(R)}}^\ast ({f^{(R)}}^\ast {\rm id}_{j_!\cF})^{\rm log}$$
in
$H_{D^+}^0(X,\cK_X).$
\end{lemma}
\begin{proof}
We consider the commutative diagram
\[\xymatrix{
{f^{(R)}}^\ast {\rm id}_{j_!\cF} \in H_{X^{(R)}}^0((X \times X)^{(R)},\cH^{(R)}) \ar[r]&
H_{X^{(R)}}^0((X \times X)^{(R)},\cH^{(R)} \otimes h^\ast R{\widetilde{g}}_\ast \La)\\
f^\ast {\rm id}_{j_!\cF} \in H_{\widetilde{X}_{U'}}^0((X \times X) \widetilde{}_{U'},\widetilde{\cH}_{U'}) \ar[r]\ar[u]^{h^\ast} &
H_{\widetilde{X}_{U'}}^0((X \times X) \widetilde{}_{U'},\widetilde{\cH}_{U'} \otimes R{\widetilde{g}}_\ast \La) \ar[u]^{h^\ast} 
}
\]
\[\xymatrix{
 & H_{X^{(R)} \setminus \widetilde{V}}^0((X \times X)^{(R)},\cH^{(R)} \otimes h^\ast R{\widetilde{g}}_\ast \La) 
 \ar[l]^{\!\!\!\!\!\!\!\!\!\!\!\!\!\!\!\!\!\!\!\!\!\!\!\!\!\!\!\!\!\!\!\!\!\!\!\!\!\!\!\!\!\!\!\!\!\!\!\!\!\!\!\!\!\!\!\!\!\!\simeq}_{\!\!\!\!\!\!\!\!\!\!\!\!\!\!\!\!\!\!\!\!\!\!\!\!\!\!\!\!\!\!\!\!\!\!\!\!\!\!\!\!\!\!\!\!\!\!\!\!\!\!\!\!\!\!\!
(\ref{map1})}\ar[r]^{{\rm Tr}^{(R)} \cdot {{\delta}^{(R)}}^\ast\!\!\!\!} &
H_{D^+}^0(X,\cK_X \otimes {\widetilde{\delta}}^\ast R{\widetilde{g}}_\ast \La) \simeq
H_{D^+}^0(X,\cK_X) \\
 & H_{\widetilde{X}_{U'} \setminus \widetilde{V}}^0((X \times X) \widetilde{}_{U'},\widetilde{\cH}_{U'} \otimes R{\widetilde{g}}_\ast \La)\ar[u]^{h^\ast}\ar[l]^{\!\!\!\!\!\!\!\!\!\!\!\!\!\!\!\!\!\!\!\!\!\!\!\!\!\!\!\!\!\!\!\!\!\!\!\!\!\!\!\!\!\!\!\!\!\!\!\!\!\!\!\!\!\!\simeq}_{\!\!\!\!\!\!\!\!\!\!\!\!\!\!\!\!\!\!\!\!\!\!\!\!\!\!\!\!\!\!\!\!\!\!\!\!\!\!\!\!\!\!\!\!\!\!\!\!\!\!\!(\ref{13'})}
 \ar[r]^{\!\!\!\!\!\!\!\!\!\!\!\!\!\widetilde{\rm Tr} \cdot {\widetilde{\delta}}^\ast} &
H_{D^+}^0(X,\cK_X \otimes {\widetilde{\delta}}^\ast R{\widetilde{g}}_\ast \La) \simeq
H_{D^+}^0(X,\cK_X)\ar[u]^{\rm id}.
}
\]
The right( resp.\ left) hand side of the equality is the image of
${f^{(R)}}^\ast {\rm id}_{j_!\cF} \in H_{X^{(R)}}^0((X \times X)^{(R)},\cH^{(R)})$ ( resp.\ $f^\ast {\rm id}_{j_!\cF} \in H_{\widetilde{X}_{U'}}^0((X \times X) \widetilde{}_{U'},\widetilde{\cH}_{U'})$) by
the composite
of the maps in the upper( resp.\ lower) lines in the commutative diagrams.
Therefore the assertion follows from the commutative diagram.
\end{proof}
We study the relationship between the elements $({f^{(R)}}^\ast {\rm id}_{j_!\cF})'$
and $({f^{(R)}}^\ast {\rm id}_{j_!\cF})^{\rm log}$ to compare the
two localizations.
\begin{lemma}\label{new}
We consider the composition 
$$H_X^0((X \times X)^{(R)},\cH^{(R)}) \longrightarrow
H_X^0((X \times X)^{(R)},\cH^{(R)} \otimes h^\ast R{\widetilde{g}}_\ast \La)
\simeq H_{D^+}^0((X \times X)^{(R)},\cH^{(R)} \otimes h^\ast R{\widetilde{g}}_\ast \La)$$
where the first map is induced by the canonical map
$\La \longrightarrow h^\ast R{\widetilde{g}}_\ast \La$ and
the second isomorphism is (\ref{map3}).
Let $(e \cup [X])^{\rm log} \in H_{D^+}^0((X \times X)^{(R)},\cH^{(R)} \otimes h^\ast R{\widetilde{g}}_\ast \La)$ 
be the image of $e \cup [X] \in H_X^0((X \times X)^{(R)},\cH^{(R)})$ by this composition. 
Then, we have an equality
$$({f^{(R)}}^\ast {\rm id}_{j_!\cF})'=({f^{(R)}}^\ast {\rm id}_{j_!\cF})^{\rm log}-(e \cup [X])^{\rm log}$$
in $H_{E^+}^0((X \times X)^{(R)},{\cH}^{(R)} \otimes h^\ast R{\widetilde{g}}_\ast \La).$
\end{lemma}
\begin{proof}
Since we have $T^{(R)} \subset E^+$, we obtain the following commutative diagram
\[\xymatrix{
({f^{(R)}}^\ast {\rm id}_{j_!\cF})' \in H_{T^{(R)}}^0((X \times X)^{(R)},{\cH}^{(R)}) \ar[d]^{\rm can.}\ar[r] &
H_{E^+}^0((X \times X)^{(R)},{\cH}^{(R)} \otimes h^\ast R{\widetilde{g}}_\ast \La)
\ar[d]_{\simeq}^{(\ref{map1})}\\
{f^{(R)}}^\ast {\rm id}_{j_!\cF}-e \cup [X] \in H_{X^{(R)}}^0((X \times X)^{(R)},{\cH}^{(R)}) \ar[r]^{\rm log\!\!\!\!\!\!\!\!\!\!\!\!\!\!\!\!\!\!\!\!} &
H_{X^{(R)}}^0((X \times X)^{(R)},{\cH}^{(R)} \otimes h^\ast R{\widetilde{g}}_\ast \La)
}
\]
where the horizontal arrows are induced by the canonical map $\La \longrightarrow
h^\ast R{\widetilde{g}}_\ast \La.$
By this commutative diagram, the isomorphism (\ref{map1})
and Corollary \ref{C1},
the assertion is proved.
\end{proof}
By Lemmas \ref{l3'} and \ref{new}, and Definition 
\ref{deg def}, it suffices to
 calculate the element $(e \cup [X])^{\rm log} \in H_{E^+}^{0}((X \times X)^{(R)},\cH^{(R)} \otimes h^\ast R{\widetilde{g}}_\ast \La)$
to prove our main theorem (Theorem \ref{theorem}).
For the calculation, we will write it
in terms of a cycle class with support on the vector bundle $E^+$
in the following Lemma \ref{l6}.

We consider the composite
\begin{align}\label{mapp}
{\small
H_X^{2d}((X \times X)^{(R)},\La(d))
\longrightarrow
H_X^{2d}((X \times X)^{(R)},h^\ast R{\widetilde{g}}_\ast \La(d))
\simeq
H_{D^+}^{2d}((X \times X)^{(R)},h^\ast R{\widetilde{g}}_\ast \La(d))}
\end{align}
where the first map is induced by the canonical map $\La \longrightarrow h^\ast R{\widetilde{g}}_\ast \La$
and the second isomorphism is (\ref{map2}).
Let $[X]^{\rm log} \in H_{D^+}^{2d}((X \times X)^{(R)},h^\ast R{\widetilde{g}}_\ast \La(d))$ denote the image
of the element
$[X] \in H_X^{2d}((X \times X)^{(R)},\La(d))$
under the map (\ref{mapp}).
By the canonical cup pairing
\begin{align}\label{pair2}
{\small\cup:\Gamma(X,{\delta^{(R)}}^\ast j^{(R)}_\ast \cH_0) \times
H_X^{2d}((X \times X)^{(R)},h^\ast R{\widetilde{g}}_\ast \La(d)) \longrightarrow
H_X^0((X \times X)^{(R)},\cH^{(R)} \otimes h^\ast R{\widetilde{g}}_\ast \La)}
\end{align} and the isomorphisms (\ref{map3}), (\ref{map2}),
we obtain an element $e \cup [X]^{\rm log}$
in $H_{D^+}^0((X \times X)^{(R)},\cH^{(R)} \otimes h^\ast R{\widetilde{g}}_\ast \La).$
By the definitions of the two pairings
(\ref{pairing1}) and (\ref{pair2}), we 
obtain an equality
\begin{equation}\label{l5}
(e \cup [X])^{\rm log}=e \cup [X]^{\rm log} 
\end{equation}
in $H_{D^+}^0((X \times X)^{(R)},\cH^{(R)} \otimes h^\ast R{\widetilde{g}}_\ast \La)$.
\begin{lemma}\label{l4}Let the notation be as above.
\\1. We consider the cartesian diagram
\[\xymatrix{
h^{-1}(X) \ar[r]\ar[d] &
(X \times X)^{(R)} \ar[d]^{h} & (V \times V) \widetilde{} \ar[l]\ar[d]^{\rm id}\\
X \ar[r]^{\!\!\!\!\!\!\!\!\!\!\!\!\!\!\!\widetilde{\delta}} &
(X \times X) \widetilde{}_{U'} & (V \times V) \widetilde{}.\ar[l]
}
\]
Let $h^!:CH_d(X) \longrightarrow CH_d(h^{-1}(X))$
be the Gysin map.
Then there exists a unique element
$[X-h^!X] \in {CH}_{d}(E^+)$
which goes to $[X]-h^![X] \in {CH}_d(h^{-1}(X))$
under the canonical map
${CH}_{d}(E^+)
\longrightarrow
{CH}_d(h^{-1}(X)).$
\\2. Further, this element satisfies an equality
$${\delta^{(R)}}^![X-h^!X]
=(-1)^{d} \cdot c_d(\Omega_{X/k}^1({\rm log}D)(R)-\Omega_{X/k}^1({\rm log}D'))_{D^+}^{X}\cap [X]$$ in $CH_{0}(D^+)$
where $c_d(.)_{D^+}^{X}\cap [X]$ is the localized Chern class introduced in
\cite[Section 3.4]{KS}.
\\3. \label{l6}We also write ${\rm cl}([X-h^!X])$ for
the image of the cycle class
${\rm cl}([X-h^!X])$ by the canonical map
$H_{E^+}^{2d}((X \times X)^{(R)},\La(d)) \longrightarrow
H_{E^+}^{2d}((X \times X)^{(R)},h^\ast R{\widetilde{g}}_\ast \La(d))$
where ${\rm cl}:{CH}_d(E^+) \longrightarrow H_{E^+}^{2d}((X \times X)^{(R)},\La(d))$ is the cycle class map.
Then we have an equality
$${\rm cl}([X-h^!X])=
[X]^{\rm log}$$ in $H_{E^+}^{2d}((X \times X)^{(R)},h^\ast R{\widetilde{g}}_\ast \La(d)).$
\end{lemma}
\begin{proof}First we prove 1,2.
We consider the commutative diagram
\[\xymatrix{
X \ar[r]^{\!\!\!\!\!\!\!\!\delta^{(R)}}\ar[d]^{id} &
(X \times X)^{(R)} \ar[d]^{h}\\
X \ar[r]^{\!\!\!\!\!\!\!\!\!\!\!\!\!\!\!\widetilde{\delta}} &
{(X \times X) \widetilde{}}_{U'}
}
\]where the horizontal arrows are the diagonal closed immersions 
 and $h:(X \times X)^{(R)} \longrightarrow (X \times X) \widetilde{}_{U'}$
is the projection.
The conormal sheaves $N_{X/(X \times X)^{(R)}}$
and $N_{X/ (X \times X) \widetilde{}_{U'}}$
are naturally identified with $\Omega_{X/k}^1({\rm log}D)(R)$
and $\Omega_{X/k}^1({\rm log}D')$ respectively, and we have
$h^{-1}(X) \setminus V=E^+.$
Therefore it is sufficient to  apply Lemma 3.4.9 in \cite{KS} to the diagram
by taking $(V \times V) \widetilde{} \subset (X \times X)^{(R)}$
as the open subscheme
$U \subset Y$ in Lemma 3.4.9 in \cite{KS}.

We prove 3.
By the following commutative diagram
\[\xymatrix{
h^\ast [X] \in H_{X^{(R)}}^{2d}((X \times X)^{(R)},\La(d)) \ar[r]^{\!\!\!\!\!\!\!\!\!\!\!\!\!\!\!\!\rm log} &
H_{X^{(R)}}^{2d}((X \times X)^{(R)},h^\ast R{\widetilde{g}}_\ast \La(d))
\ni (h^\ast [X])^{\rm log}\\
[X] \in H_X^{2d}((X \times X) \widetilde{}_{U'},\La(d)) \ar[r]\ar[u]^{h^\ast} &
H_X^{2d}((X \times X) \widetilde{}_{U'},R{\widetilde{g}}_\ast \La(d))=0, \ar[u]^{h^\ast}}
\]
we acquire $(h^\ast [X])^{\rm log}=0$.
We consider the commutative diagram
\[\xymatrix{
{\rm cl}([X-h^!X]) \in H_{E^+}^{2d}((X \times X)^{(R)},\La(d)) \ar[r]\ar[d]^{(1)} &
H_{E^+}^{2d}((X \times X)^{(R)},h^\ast R{\widetilde{g}}_\ast \La(d))  \ar[d]_{\simeq}^{(\ref{map4})}\\
[X]-h^\ast [X] \in H_{X^{(R)}}^{2d}((X \times X)^{(R)},\La(d)) \ar[r]^{\rm log\!\!\!\!\!\!\!\!} &
H_{X^{(R)}}^{2d}((X \times X)^{(R)},h^\ast R{\widetilde{g}}_\ast \La(d)).
}
\]
The element
${\rm cl}([X-h^!X]) \in H_{E^+}^{2d}((X \times X)^{(R)},\La(d))$ goes to $[X]-h^\ast [X] \in H_{X^{(R)}}^{2d}((X \times X)^{(R)},\La(d))$ under the map (1) by Lemma \ref{l4}.
Hence the assertion follows from this commutative diagram and $(h^\ast [X])^{\rm log}=0$.
\end{proof}

We are ready to prove the main result in this paper.
\begin{theorem}\label{theorem}
Let $X$ be a smooth scheme of dimension $d$
over a perfect field $k$
and $U \subset X$ be the complement 
of a divisor with simple normal crossings.
Let $\cF$ be a smooth sheaf of rank 1, and $T$ and $D^+$
the nonclean locus and the wild locus of $\cF$ respectively.
Then we have an equality
$$C_{T}^{{\rm ncl},0}(j_! \cF)+
(-1)^{d} \cdot c_d(\Omega_{X/k}^1({\rm log}D)(R)-\Omega_{X/k}^1({\rm log}D'))_{D^+}^{X}\cap [X]
=C_{D^+}^{{\rm log},0}(j_! \cF)$$
in $H_{D^+}^0(X,\cK_X).$
\end{theorem}
\begin{proof}
By Lemmas \ref{l3'} and \ref{new}, and the commutative diagram
\[\xymatrix{
({f^{(R)}}^\ast {\rm id}_{j_!\cF})' \in H_{T^{(R)}}^0((X \times X)^{(R)},{\cH}^{(R)}) \ar[r]\ar[d]^{{\rm Tr}^{(R)} \cdot {\delta^{(R)}}^\ast}
& H_{E^+}^0((X \times X)^{(R)},{\cH}^{(R)} \otimes h^\ast R{\widetilde{g}}_\ast \La)\ar[d]^{{\rm Tr}^{(R)} \cdot {\delta^{(R)}}^\ast} \\
\!\!\!\!\!\!\!\!\!\!\!\!\!\!\!\!\!\!\!\!\!\!\!\!\!\!\!\!\!\!\!\!\!\!\!\!\!\!\!\!\!\!\!\!\!C_{T}^{{\rm ncl},0}(j_! \cF) \in H_T^0(X,\cK_X) \ar[r]^{\!\!\!\!\!\!\!\!\!\!\!\!\!\!\!\!\!\!\!\!\!\!\!\!\!\!\!\!\!\!\!\!\!\!\!\!\!\!\!\!\!\!\!\!\!\!\!\!\!\!\rm can.} &
H_{D^+}^0(X,\cK_X) \simeq H_{D^+}^0(X,\cK_X \otimes {\widetilde{\delta}}^\ast R{\widetilde{g}}_\ast \La),
}
\]
we obtain an equality 
$$C_{T}^{{\rm ncl},0}(j_! \cF)=C_{D^+}^{{\rm log},0}(j_! \cF)-
{{\rm Tr}^{(R)}} \cdot {\delta^{(R)}}^\ast ((e \cup [X])^{\rm log})$$
in $H_{D^+}^0(X,\cK_X).$
Since we have  
${{\rm Tr}^{(R)}} {\delta^{(R)}}^\ast ({e} \cup [X]^{\rm log})
={{\rm Tr}^{(R)}}({e}) \cup {\delta^{(R)}}^\ast
[X]^{\rm log}={\delta^{(R)}}^\ast
[X]^{\rm log}$ in $H_{D^+}^0(X,\cK_X),$ 
we acquire equalities
$${{\rm Tr}^{(R)}} \cdot {\delta^{(R)}}^\ast ((e \cup [X])^{\rm log})
={{\rm Tr}^{(R)}} \cdot {\delta^{(R)}}^\ast (e \cup [X]^{\rm log})
={\delta^{(R)}}^\ast {\rm cl}([X-h^!X])$$
$$=(-1)^{d} \cdot c_d(\Omega_{X/k}^1({\rm log}D)(R)-\Omega_{X/k}^1({\rm log}D'))_{D^+}^{X}\cap [X]$$
by the equality (\ref{l5}), and Lemma \ref{l4}.
Therefore the assertion is proved.
\end{proof}
We show that the logarithmic localized characteristic class
refines the characteristic class.
The localized Chern class
$(-1)^{d} \cdot c_d(\Omega_{X/k}^1({\rm log}D)(R)-\Omega_{X/k}^1({\rm log}D'))_{D^+}^{X}\cap [X] \in H_{D^+}^0(X,\cK_X)$
goes to the difference
$(-1)^{d} \cdot c_d(\Omega_{X/k}^1({\rm log}D)(R)) \cap [X]-
(-1)^{d} \cdot c_d(\Omega_{X/k}^1({\rm log}D')) \cap [X] \in H^0(X,\cK_X)$
by the canonical map $H_{D^+}^0(X,\cK_X)
\longrightarrow H^0(X,\cK_X).$
(See \cite[Section 3.4]{KS}.)

Let $(\Delta_X,\Delta_X) \in H^0(X,\cK_X)$ denote the self-intersection product.\ Since the conormal sheaf $N_{X/(X \times X)^{(R)}}$
( resp.\ $N_{X/(X \times X) \widetilde{}_{U'}}$)
is isomorphic to $\Omega_{X/k}^1({\rm log}D)(R),$
( resp.\ $\Omega_{X/k}^1({\rm log}D')$,)
we have an equality
$(-1)^{d} \cdot c_d(\Omega_{X/k}^1({\rm log}D)(R)) \cap [X]
=(\Delta_X,\Delta_X)_{(X \times X)^{(R)}}.$
( resp.\ $(-1)^{d} \cdot c_d(\Omega_{X/k}^1({\rm log}D')) \cap [X]=
(\Delta_X,\Delta_X)_{(X \times X) \widetilde{}_{U'}}.$)
\begin{corollary}\label{coro}
Let the notation be as in Theorem \ref{theorem}.
Let $C(.)$ denote the characteristic class defined in \cite[Definition 2.1.1]{AS}.
Then, the image of the logarithmic localized characteristic class
$C_{D^+}^{\rm log,0}(j_!\cF) \in H_{D^+}^0(X,\cK_X)$
under the canonical map $H_{D^+}^0(X,\cK_X) \longrightarrow H^0(X,\cK_X)$
is the difference $C(j_!\cF)-C({j_{U'}}_!\La_{U'}).$
\end{corollary}
\begin{proof}
By Definition \ref{deg def},
the nonclean localized characteristic class $C_{T}^{{\rm ncl},0}(j_! \cF)$ 
 goes to
the difference
$C(j_!\cF)-(\Delta_X,\Delta_X)_{(X \times X)^{(R)}} \in H^0(X,\cK_X)$
by the canonical map
$H_T^0(X,\cK_X) \longrightarrow H^0(X,\cK_X).$
The localized Chern class
$(-1)^{d} \cdot c_d(\Omega_{X/k}^1({\rm log}D)(R)-\Omega_{X/k}^1({\rm log}D'))_{D^+}^{X}\cap [X]$
goes to the difference
$$(\Delta_X,\Delta_X)_{(X \times X)^{(R)}}\\-(\Delta_X,\Delta_X)_{(X \times X) \widetilde{}_{U'}}$$
under the canonical map
$H_{D^+}^0(X,\cK_X) \longrightarrow H^0(X,\cK_X).$
We have an equality $C({j_{U'}}_!\La_{U'})=(\Delta_X,\Delta_X)_{(X \times X) \widetilde{}_{U'}}$ in $H^0(X,\cK_X)$ by \cite[Corollary 2.2.5.1]{AS}.
Hence the assertion follows from Theorem \ref{theorem}.
\end{proof}

Let the notation be as above. We assume that $\cF$ is a smooth
$\La$-sheaf of rank 1 which is clean with respect to the boundary $D$. 
Let $c_{\cF} \in CH_0(D^+)$
be the Kato 0-cycle class of $\cF$. (See \cite{K2} or \cite[Definition 4.2.1]{AS}.)
We have an equality 
$c_{\cF}=-(-1)^d{c_d}(\Omega^1_{X/k}({\rm log} D)(R)-\Omega^1_{X/k}({\rm log} D) 
)_{D^+}^{X}\cap [X]$.
(See loc. cit.) 
\begin{corollary}\label{corollary}
Let $X$ be a smooth scheme of dimension $d$
over a perfect field $k$
and $U \subset X$ be the complement 
of a divisor with simple normal crossings.
Let $\cF$ be a smooth sheaf of rank 1 and $D^+$ the wild
locus of $\cF.$
If $\cF$ is clean with respect to the boundary, we have
$$-c_{\cF}=C_{D^+}^{{\rm log},0}(j_! \cF)-C_{D^+}^{{\rm log},0}(j_! \La_U)$$
in $H_{D^+}^0(X,\cK_X).$
\end{corollary}
\begin{proof}
Since we assume that $\cF$ is clean with respect to the boundary $D$,
we have $C_T^{\rm ncl,0}(j_! \cF)=0.$
Thereby
the equalities
$(-1)^{d}c_d(\Omega_{X/k}^1({\rm log}D)(R)-\Omega_{X/k}^1({\rm log}D'))_{D^+}^{X}\cap [X]=C_{D^+}^{{\rm log},0}(j_! \cF)$ and
$(-1)^{d}c_d(\Omega_{X/k}^1({\rm log}D)-\Omega_{X/k}^1({\rm log}D'))_{D^+}^{X}\cap [X]=C_{D^+}^{{\rm log},0}(j_! \La_U)$ hold in $H_{D^+}^0(X,\cK_X)$
by Theorem \ref{theorem}.
Hence we obtain
$$C_{D^+}^{{\rm log},0}(j_! \cF)-C_{D^+}^{{\rm log},0}(j_! \La_U)
=(-1)^d{c_d}(\Omega^1_{X/k}({\rm log} D)(R)\\-\Omega^1_{X/k}({\rm log} D) 
)_{D^+}^{X}\cap [X]=-c_{\cF}.$$
\end{proof}

\end{document}